\documentclass[letterpaper, 10 pt, conference]{ieeeconf}  
\IEEEoverridecommandlockouts        
\usepackage{graphicx} 
\usepackage{epstopdf} 
\usepackage{psfrag} 
\usepackage{amsmath} 
\usepackage{amssymb}  
\usepackage{amsthm,thmtools,dsfont}
\usepackage{color}
\usepackage{tabu}
\usepackage{hyperref,cleveref}
\let\labelindent\relax
\usepackage{enumitem} \setlist[itemize,1]{leftmargin=\dimexpr 22pt}
\usepackage{cite}
\usepackage{microtype,enumerate,url}
\usepackage{BOONDOX-ds}
\usepackage{tikz} 
\usepackage[fulladjust]{marginnote}
\usepackage{pgfplotstable} 
\usepackage{xstring} 
\usepackage{booktabs} 
\usepackage{colortbl}
\usepackage{circuitikz}
\usetikzlibrary{arrows,patterns}
\usepackage{caption}
\usepackage{stfloats}%

\usepackage{fontawesome}
\newcommand{\newcheckmark}{\textrm{\faCheck}}
\newcommand{\newcrossmark}{\textrm{\faTimes}}

\definecolor{Gray}{gray}{0.9}

\usepackage{algorithm}
\usepackage{algorithmic}

\renewcommand{\theenumi}{(\roman{enumi})}

\makeatletter
\def\@IEEEsectpunct{.\ \,}
\def\paragraph{\@startsection{paragraph}{4}{\z@}{1.5ex plus 1.5ex minus 0.5ex}%
{0ex}{\normalfont\normalsize\itshape}}
\makeatother

\renewcommand*\descriptionlabel[1]{\hspace\labelsep\normalfont#1}

\declaretheorem[style=definition]{theorem}

\declaretheorem[style=definition]{lemma}
\declaretheorem[style=definition,qed=$\vartriangle$]{remark}

\declaretheorem[style=definition,numbered=no]{standing assumption}

\declaretheorem[style=definition]{definition}

\renewcommand\thmcontinues[1]{continued}

\newcommand {\nn}{\nonumber}
\newcommand{\beq}{\begin{equation}}
\newcommand{\eeq}{\end{equation}}
\newcommand {\bseq}{\begin{subequations}}
\newcommand {\eseq}{\end{subequations}}
\newcommand {\bma}{\left[}
\newcommand {\ema}{\right]}

\newcommand {\N}{\mathbb{N}} 	
\newcommand {\Zge}{\mathbb{Z}_{+}} 	
\newcommand {\R}{\mathbb{R}} 	
\newcommand {\Rge}{\mathbb{R}_{+}} 	

\newcommand {\T}{\mathbb{T}} 	
\newcommand {\W}{\mathbb{W}} 	
\newcommand {\B}{\mathcal{B}} 	
\newcommand {\C}{\mathcal{C}} 	
\newcommand {\K}{\mathcal{K}} 	
\newcommand {\A}{\mathcal{A}} 	
\newcommand {\M}{\mathcal{M}} 	
\renewcommand{\L}{\mathcal{L}} 	

\newcommand{\Span}{\mathbf{span}}
\newcommand{\aff}{\mathbf{aff}} 
\newcommand{\cone}{\mathbf{cone}} 
\newcommand{\ccone}{\mathbf{cvxcone}} 
\newcommand{\conv}{\mathbf{conv}} 
\newcommand{\closure}{\mathbf{cl}} 
\newcommand{\rank}{\mathbf{rank}} 
\newcommand{\transpose}{\mathsf{T}} 

\ctikzset{
    resistors/scale=0.7,
    capacitors/scale=0.7,
    diodes/scale=0.6,
    inductors/coils=6,
    blocks/scale=0.6
}

\title{\Large {\bfseries   Data-driven representations of  conical, convex, and affine  behaviors}}

\author{Alberto Padoan, Florian D\"orfler, John Lygeros\thanks{
		A.  Padoan, F. D\"orfler, J. Lygeros are with 
		the Department of Information Technology and Electrical Engineering at		
		ETH Z\"urich,  Z\"urich, Switzerland {\tt\footnotesize  \{apadoan, dorfler, lygeros\}@control.ee.ethz.ch}. 
		Research supported by the Swiss National Science Foundation under the NCCR Automation.
  }
}

\date{\small\today} 

\begin{document}

\maketitle
\thispagestyle{empty}
\pagestyle{empty}

\begin{abstract}   
\noindent  The paper studies  conical, convex, and affine    models in the framework of behavioral systems theory.  
We  investigate basic properties of  such behaviors  
 and address  the problem of constructing models from measured data. 
We prove that  closed, shift-invariant,  conical, convex, and affine   models have the intersection  property,  
thereby enabling the definition of most powerful unfalsified models based on infinite-horizon measurements. 
We then provide  necessary and sufficient  conditions for representing  conical, convex, and affine finite-horizon  behaviors using raw data matrices, expressing persistence of excitation requirements in terms of non-negative rank conditions. The applicability of our results is demonstrated by a numerical example arising  in population ecology. 
\end{abstract}

\section{Introduction}
 
The behavioral approach to systems modeling is a milestone of systems and control theory.
Introduced in the landmark papers~\cite{willems1986timea,willems1986timeb,willems1987timec},  
it  has made substantial contributions to the field, decoupling system models from their representations, resolving long-standing  modeling challenges  related to ports and input-output interconnections~\cite{willems2010ports}, and recently enabling a new wave of data-driven control algorithms via the fundamental lemma~\cite{willems2005note,markovsky2021behavioral}.  

At the heart of behavioral systems theory  is the core principle that a system is characterized by its behavior, i.e., the set of all possible trajectories of the system.  A system is linear if the corresponding behavior is a subspace.  Over time, the   focus of behavioral systems theory has been primarily directed to  finite-dimensional, linear, time-invariant (LTI) systems~\cite{willems1986timea,willems1986timeb,willems1987timec}, whose behaviors are closed, shift-invariant, subspaces. However,  behavioral systems theory remains   largely unexplored for  other   classes of systems.

The  objective of this paper is to challenge the standard assumption of linearity and examine the situation where the behavior of a system is a conical, convex, or affine set.  
Our investigation centers around the idea that   such behaviors  provide a powerful and tractable modeling tool.  Conical and convex behaviors capture global information from sparse or limited data, often leading to simple representations. We illustrate this principle by considering the problem of building  conical, convex, and affine  models from measured data,  both over finite and   infinite time  horizons.

\textbf{Contributions}:   The main contributions of the paper can be summarized as follows.  
(i) We  define  systems with  conical, convex, and affine  behaviors and show that they possess a number of desirable properties. (ii) We study the construction of  conical, convex, and affine  models from given measurements using the notion of \textit{most powerful unfalsified model}~\cite{willems1986timeb}. We prove that   
discrete-time  models with closed, shift-invariant, behaviors that are either conical  (e.g., positive),  convex, or affine have the intersection property, thereby enabling the definition of a most powerful unfalsified model~\cite{willems1986timeb}. (iii) We provide necessary and sufficient conditions echoing Willems' fundamental lemma~\cite{willems2005note}  for representing finite-horizon behaviors of time-invariant, positive, linear and affine state-space systems using a Hankel matrix constructed from measured data, expressing persistence of excitation requirements in terms of non-negative  low  rank conditions.  (iv) Finally, we illustrate with a numerical case study that prior information about the conicity  and convexity of a system facilitates the construction of models from data.

\textbf{Related work}: 
The present paper builds on the notions of conical and convex systems introduced in~\cite{dam1997unilaterally}, 
expanding their scope to  conical  behaviors that are not necessarily polyhedral. 
Willems' fundamental lemma~\cite{willems2005note} gives conditions for the image of a Hankel matrix constructed from a trajectory of the system to be a data-driven representation of the 
behavior of a discrete-time LTI system  over a finite time horizon.  
An overview of different data-driven representations  of  finite-dimensional LTI systems can be found in~\cite{markovsky2021behavioral}.  Data-driven   representations of affine systems have been discussed, e.g., and~\cite{martinelli2022data}. A data-driven approach to the stabilization of positive linear state-space systems has been  studied  in~\cite{shafai2022data}. A geometric view on data-driven behavioral systems theory has been explored in~\cite{padoan2022behavioral}.

\textbf{Paper organization}:  The rest of the paper is organized as follows.  
Section~\ref{sec:behaviors} focuses on   conical, convex, and affine  behaviors,
covering their definition and elementary properties.
Section~\ref{sec:main_results} discusses the construction of conical, convex, and affine models from given measurements. 
Section~\ref{sec:example} illustrates the theory with a numerical example. 
Section~\ref{sec:conclusion} provides a summary of our main results and an outlook to future research directions. 
The proofs of our main results are deferred to the appendix.

\textbf{Notation}:  
The sets of positive and non-negative integers are denoted by $\N$ and $\Zge$, respectively.   
The sets of real and non-negative real numbers are denoted by $\R$ and $\Rge$, respectively.    
For ${T\in\N},$ the set of integers $\{1, 2, \dots , T\}$ is denoted by $\mathbf{T}$.  
A map $f$ from $X$ to $Y$ is denoted by ${f:X \to Y}$; $(Y)^{X}$ denotes the set of all such maps. 
The restriction of $f:X \to Y$ to a set $X^{\prime}$, with $X^{\prime}\cap X \not = \emptyset$, is denoted by 
$f|_{X^{\prime}}$ and is defined by $f|_{X^{\prime}}(x) = f(x)$ for ${x \in X^{\prime} \cap X}$.
If $\mathcal{F} \subseteq (Y)^{X}$, then  $\mathcal{F}|_{X^{\prime}}$ denotes ${\{ f|_{X^{\prime}} \, : \, f \in \mathcal{F}\}}$. The vector whose entries are all one is denoted by $\mathbb{1}$. 
The rank of a matrix ${M\in\R^{p\times m}}$  is denoted by $\rank\, M $.

\section{Conical, convex, and affine behaviors} \label{sec:behaviors}

In the language of behavioral systems theory~\cite{willems1986timea,willems1986timeb,willems1987timec},
a \textit{dynamical system}  (or, briefly, a \textit{system})   is a triple $\Sigma=(\T,\W,\B),$ where $\T$ is the \textit{time set}, $\W$ is the \textit{signal set}, and $\B \subseteq (\W)^{\T}$ is the \textit{behavior} of the dynamical system.  
The time axis $\T$ is (usually) $\Zge$ for \textit{discrete-time} systems or (sometimes) $\Rge$ for \textit{continuous-time} systems. The signal set $\W$ is assumed to be a linear space, unless otherwise stated. 
By a convenient abuse of terminology, we often identify a system with the corresponding behavior. 

A \textit{trajectory} is any element ${w\in(\W)^{\T}}$.  For discrete-time systems, a trajectory of length ${T\in\N}$
is any element ${w\in(\W)^{\Zge}|_{[0,T-1]}}$. By a convenient abuse of notation, we often identify trajectories of length ${T\in\N}$ with the corresponding vector  ${w=(w(0),\ldots,w(T-1)) \in \W^{T}}$.

Given ${\tau\in\T}$ and a trajectory ${w\in(\W)^{\T}}$, the \textit{(backwards) $\tau$-shift} is defined as $(\sigma^\tau w)(\tau) = w(t+\tau)$ for all ${t\in\T}$. A discrete-time system $\Sigma$ is \textit{shift-invariant} if the corresponding behavior $\B$ is such that ${\sigma^t(\B) \subseteq \B}$, for all ${t \in \Zge},$  and \textit{complete} if the condition ${w\in\B}$  holds if and only if $\,w|_{[t_1,t_2]} \in\B|_{[t_1,t_2]}$,  for all ${t_1, t_2 \in \Zge}$, with  ${t_1\le t_2}$.  

\subsection{Conical, convex, affine, and linear systems}

A subset $K$  of $\W$ is a \textit{cone} if ${\alpha \in\Rge}$ and ${x\in K}$ imply $\alpha x\in K$~\cite[p.1]{zalinescu2002convex}.
A subset ${C}$ of $\W$   is \textit{convex} if ${\alpha \in [0,1]}$ and ${x,y \in C}$ imply ${\alpha x + (1-\alpha) y \in C}$~\cite[p.1]{zalinescu2002convex}.  A subset ${  A   \subseteq \W}$ is a \textit{affine} if ${\alpha \in \R}$ and ${x,y \in A}$  imply ${\alpha x + (1-\alpha) y \in  A  }$~\cite[p.1]{zalinescu2002convex}.  
The conical, convex, convex conical, affine, and linear hull of a subset $S$ of $\W$  are defined as~\cite[p.2]{zalinescu2002convex}
\begin{align*} 
\cone \, S &= \bigcap \, \{ K \subseteq \W \,:\,  S \subseteq K, \ K \text{ cone}  \}, \\ 
\conv \, S &= \bigcap \, \{ C \subseteq \W \,:\,  S \subseteq C, \ C \text{ convex set}  \}, \\ 
\ccone \, S &= \bigcap \, \{ V \subseteq \W \,:\,  S \subseteq V, \  V \text{ convex cone} \}, \\ 
\aff \, S &= \bigcap \, \{ A \subseteq \W \,:\,  S \subseteq A, \ A \text{ affine set}  \}, \\ 
\Span  \, S &= \bigcap \,  \{ L \subseteq \W \,:\,  S \subseteq L, \  L \text{ linear space}  \} , 
\end{align*}  
respectively.

\begin{definition}
A dynamical system $\Sigma=(\T,\W,\B)$ is \textit{conical (convex, affine, linear)} if $\B$ is a cone (convex, affine, linear).
\end{definition}

\noindent
Convex and conical behaviors have been introduced in~\cite{dam1997unilaterally} with the aim of modeling constrained mechanical systems. However,  they can   model a much wider class of systems.

\subsection{Examples}

\subsubsection{Linear systems}
Linear systems are obviously conical, convex, and affine.  This includes  
infinite-dimensional linear systems  (e.g., ${y(t) = u(t-1)}$ in continuous-time)  as well as finite-dimensional LTI systems, such as any discrete-time, LTI, state-space system described by the equations
\beq \label{eq:state-space-linear}  
\sigma x = Ax + Bu, \quad  y=Cx+Du, 
\eeq
with ${x(t)\in\R^n}$, ${u(t)\in\R^m}$, and ${y(t)\in\R^p}$, whose behavior 
$\B = \left\{\,
		(u,y)  \in (\R^{p+m})^{\Zge} 
		\, : \,	\exists \,  x \in (\R^{n})^{\Zge}  \text{ s.t.} ~\eqref{eq:state-space-linear}  \text{ holds}
	  \,\right\}$
is always a complete, shift-invariant subspace~\cite{willems1986timea}.

\subsubsection{Positive systems}
Positive systems are  broadly   defined by behaviors which leave a  cone   invariant~\cite{farina2000positive}.  
For example,  the system~\eqref{eq:state-space-linear} 
is \textit{positive} if the non-negative orthant is invariant under  its  dynamics~\eqref{eq:state-space-linear}, i.e., if ${x(0) \in \Rge^n}$ and ${u(t) \in \Rge^m}$  imply ${x(t) \in \Rge^n}$   and ${y(t) \in \Rge^p}$   for all ${t\in\Zge}$. 
Positivity of system~\eqref{eq:state-space-linear} ensures that  its   restriction to the non-negative orthant  is well-defined;  conicity and convexity follow directly from the properties of the non-negative orthant.

\subsubsection{Sector conditions}
Sector conditions and, more generally, dissipation inequalities model nonlinear, uncertain and time-varying behaviors of as subsets of a cone~\cite{khalil1996nonlinear}. A function $\varphi:\R\to\R$   satisfies   the \textit{sector condition} ${\varphi \in  [\alpha, \beta]}$ if
\beq \label{eq:sector-condition} 
(\varphi (u) - \alpha  u) ( \varphi (u) - \beta u) \leq 0, \  \forall \, u \in \R,
\eeq
with ${-\infty \le \alpha < \beta \le \infty}$.
The sector condition ${\varphi \in [\alpha, \beta]}$ defines a memoryless system whose  behavior 
$\B = \{(u,y)\in (\R^{p+m})^{\T}:(y - \alpha  u) ( y - \beta u) \le 0 \}$
is a (quadratic) cone.

\subsection{Properties}

We now turn our attention  to conical, convex, and affine behaviors. We begin with a list of basic properties directly inherited  from  the  underlying  structure.  

\begin{lemma}[Algebraic and topological properties]\label{lemma:algebraic_properties}
The following are conical (convex, affine, linear) behaviors.
\begin{itemize}
\item[(i)] The intersection of any family of conical (convex, affine, linear) behaviors.
\item[(ii)] The image and preimage of a conical (convex, affine, linear) behavior under a linear mapping. 
\item[(iii)] The Cartesian product of conical (convex, affine, linear) behaviors. 
\item[(iv)] The sum of conical (convex, affine, linear) behaviors.
\item[(v)]  If $(\W)^{\T}$ is a topological space,   the closure of conical (convex, affine, linear) behaviors.
\end{itemize}
\end{lemma} 

The structure of a behavior can be also characterized in terms of the elements it contains.  Recall   that a \textit{ray} is a cone of the form 
${K=\{\lambda w \, :  \lambda \in \Rge \}},$
with ${w \not = 0}$~\cite[p.15]{rockafellar1970convex}.  A \textit{linear combination} of (finitely many) elements ${w_1, \ldots,  w_n \in \W}$ is an expression of the form 
${g_1 w_1+ \ldots +  g_n w_n}$,
with ${g = (g_1, \ldots,  g_n) \in \R^n}$. \textit{Conical}, \textit{affine}, and \textit{convex combinations}
are expressions of the same form 
such that ${g \in \Rge^n},$ ${\mathbb{1}^{\transpose}g = 1}$, and both ${\mathbb{1}^{\transpose}g = 1}$ and ${g \in \Rge^n}$, respectively.  

\begin{lemma}[Combinations]\label{lemma:combinations}
The following statements hold.
\begin{itemize}
\item[(i)] A behavior is conical if and only if it contains its rays.
\item[(ii)] A behavior is linear (affine, convex, or convex conical) if and only if it contains all linear (affine, convex, or conical) combinations of its elements.
\end{itemize}
\end{lemma}

\noindent 
Time shifts and restrictions over finite time horizons  play a key role in behavioral systems theory~\cite{willems1986timea,willems1986timeb,willems1987timec}. The next result links these operations with the properties discussed above.

\begin{lemma} \label{lemma:restrictions_and_time_shifts} (Restrictions and time shifts) 
Let $\Sigma = (\T, \W, \B)$ be a conical (convex, affine, linear) system. Then, for all $t \in \T$ and all $t_1, t_2 \in \T$ such that $t_1 \leq t_2$, both $\sigma^{t}(\B)$ and $\B|_{[t_1, t_2]}$ are a cone (convex set, affine set, linear space). 
\end{lemma}

\section{From time series to discrete-time \\ conical, convex, affine, and linear systems} \label{sec:main_results}

In behavioral systems theory, a \textit{model} of a system $\Sigma$ is simply a subset $\B$ of the set of all possible trajectories $\W^{\T}$~\cite{willems1986timeb}. 
In the same spirit, a \textit{model class} $\M$ is any family of subsets of $\W^{\T}$~\cite{willems1986timeb}. In this context, the \textit{most powerful unfalsified model}   describes the ``element in a model class which explains a given set of observations and as little else as possible''~\cite[p.1]{willems1986timeb}.   The most powerful unfalsified model   plays a fundamental role in the solution of the system identification problem~\cite{markovsky2021behavioral}.  

\subsection{Most\,powerful\,unfalsified\,conical,\,convex,\,and\,affine\,models} 

Given a model class ${\M}$ and a trajectory ${w_d\in\W^{\T}},$ 
 the corresponding identified  model ${\B_{\textup{mpum}}(w_d)}$ is the \textit{most powerful unfalsified model in the model class $\mathcal{M}$ based on the measurement $w_d$} if~\cite[Definition 4]{willems1986timeb}
\begin{itemize}
\item[(i)] ${\B_{\textup{mpum}}(w_d) \in \M}$,
\item[(ii)] ${w_d \in \B_{\textup{mpum}}(w_d)}$, and
\item[(iii)] ${\B \in \M}$ and ${w_d \in \B}$ imply ${\B_{\textup{mpum}}(w_d) \subseteq \B}$.
\end{itemize} 
The most powerful unfalsified model need not exist, but if it exists it is unique~\cite{willems1986timeb}.
The most powerful unfalsified model in the model class $\L^q$ of discrete-time, finite-dimensional, LTI systems over ${\W=\R^q}$ based on the measurement ${w_d\in(\R^q)^{\Zge}}$ is~\cite{willems1986timeb} 
\beq \label{eq:MPUM_linear}
\B_{\text{mpum}}(w_d) = \closure\left( \Span\{w_d, \sigma w_d, \sigma^2 w_d, \ldots \}\right),
\eeq 
where $\closure\,(S)$ is used to denote the closure of the set $S$.

We now establish an analogous result for the model classes $\K^q$, $\C^q$, and $\A^q$ of all
discrete-time, closed, shift-invariant  models that are conical, convex, and affine over ${\W=\R^q}$, respectively.   
The distinction between the most powerful unfalsified models   in the model class $\L^q$ and 
those in $\K^q$, $\C^q$, and $\A^q$ is disarmingly simple:
 it suffices  to replace $\Span$ with $\cone$, $\conv$, and  $\aff$, respectively.  
The key technical tool behind this result is the next lemma.

\begin{lemma}[Intersection property] \label{lemma:intersection_property}
The model classes $\K^q$, $\C^q$, and $\A^q$ have the intersection  property,  i.e.,  the intersection of   any   collection of elements of $\mathcal{M}$ is  itself  an element of   $\mathcal{M}$. 
\end{lemma}

We are now ready to state the main result of this section.

\begin{theorem}[Most powerful unfalsified models] \label{thm:MPUM}
The  most powerful unfalsified models   in the model classes $\K^q$, $\C^q$, and $\A^q$ based on the measurement ${w_d\in(\R^q)^{\Zge}}$ are
\begin{align}
\B_{\text{mpum}}(w_d) &= \closure\left( \cone\{w_d, \sigma w_d, \sigma^2 w_d, \ldots \}\right),\label{eq:MPUM_cone} \\
\B_{\text{mpum}}(w_d) &= \closure\left( \conv\{w_d, \sigma w_d, \sigma^2 w_d, \ldots \}\right),\label{eq:MPUM_conv}  \\
\B_{\text{mpum}}(w_d) &= \closure\left( \aff\{w_d, \sigma w_d, \sigma^2 w_d, \ldots \}\right). \label{eq:MPUM_aff} 
\end{align}
\end{theorem}

\begin{remark}[Complete and closed behaviors]
The relationship between completeness and closedness for conical, convex, and affine systems has not been explored in our discussion. While completeness and closedness are equivalent for certain classes of discrete-time systems, including linear systems~\cite[p.567]{willems1986timea} and finite-polyhedral\footnote{A behavior $\B\subseteq (\R^q)^{\Zge}$ is \textit{finite-polyhedral} if  its restriction  to every finite-length interval of $\Zge$ is a polyhedral cone~\cite[p.44]{dam1997unilaterally}.} systems~\cite[p.44]{dam1997unilaterally}, it is unclear for what other classes of systems this relationship holds. 
\end{remark}

In the next section, we investigate finite-horizon, time-invariant, linear and affine behaviors with the additional structural constraints enforced by positivity.

\subsection{On conical, convex, and affine models over finite-horizons}

The cornerstone of several recent developments in the area of data driven control is a basic  principle: 
given a discrete-time LTI system ${\B \in \mathcal{L}^{q}}$ and a trajectory of the system ${w_d \in \R^{qT}}$, the restricted behavior $\B|_{[0,L-1]}$ can be represented by a raw data matrix, provided the time horizon ${L}$ is long enough and the data $w_d$ are sufficiently informative.

We now  revisit  a version of this principle from~\cite{markovsky2022identifiability}.   
Recall that the \textit{order} of an LTI system ${\B\in\L^q}$ with  (minimal)  state-space representation~\eqref{eq:state-space-linear}    is    the smallest ${n\in\N}$ among all representations~\eqref{eq:state-space-linear} and the \textit{lag}     is   the smallest ${\ell\in\N}$  such  that (in a minimal representation)  
the \textit{observability matrix} 
\beq  \label{eq:observability}  
\mathsf{O}_\ell
=
\bma
\begin{array}{cc}
C \\
CA \\
\vdots  \\
CA^{\ell -1} 
\end{array} 
\ema 
\eeq 
is  full rank. 
The \textit{Hankel matrix} of depth ${L\in\mathbf{T}}$ associated with the trajectory ${w_d \in \R^{qT}}$ is\,defined\,as 
\beq \label{eq:Hankel} 
\! H_{L}(w_d) \! = \!
\scalebox{0.9}{$
\bma  
\begin{array}{ccccc}
w_d(0) & w_d(1)  & \cdots &  w_d(T-L)   \\
w_d(1) & w_d(2)  & \cdots &  w_d(T-L+1)   \\
\vdots  & \vdots  & \ddots & \vdots  \\
w_d(L-1) & w_d(L)  & \cdots  & w_d(T-1)
\end{array}
\ema
$} . \! \!
\eeq  
With abuse of notation, we denote the linear, affine, convex, and convex conical hull of the columns of a matrix ${M\in\R^{p\times m}}$ by $\Span(M)$, $\aff(M)$,  $\conv(M)$, and $\ccone(M)$, respectively.

\begin{lemma}~\cite[Corollary 21]{markovsky2022identifiability} \label{lemma:fundamental_generalized}  
Let ${\B\in\L^q}$ be an LTI  system of order $n$ and lag $\ell$ with state-space representation~\eqref{eq:state-space-linear}. Let ${w_d  = (u_d,y_d)   \in \R^{qT}}$ be a trajectory of the system  and let ${L\in\mathbf{T}}$, with ${L > \ell}$. 
Then
\beq  \label{eq:BL} 
\B|_{[0,L-1]} = \Span \, H_L(w_d) 
\eeq
if and only if
\beq  \label{eq:generalized_persistence_of_excitation} 
\rank \, H_L(w_d)   =  mL+ n.
\eeq
\end{lemma}

\noindent 
The rank condition ${\rank \, H_L(w_d) =  mL+ n}$ is referred to as the \emph{generalized persistence of excitation} in~\cite{markovsky2022identifiability}. 
By the fundamental lemma~\cite{willems2005note},   the condition is guaranteed to hold for controllable systems   
if the input is  persistently exciting of order $n+L$, in which case the rank condition
\beq \nn
\rank \, 
\scalebox{1}{$
\bma  
\begin{array}{c}
H_L(u_d) \\
H_1(x_d)
\end{array}
\ema = mL+n , $}
\eeq
holds, where $x_d$ is the state trajectory of~\eqref{eq:state-space-linear} corresponding to the input trajectory $u_d$. 
An overview of different variations of this principle can be found in the recent  survey~\cite{markovsky2021behavioral}.

A conceptually similar result  has been presented in~\cite{martinelli2022data} for   \textit{affine},  \textit{time-invariant} (ATI)  systems described by the state-space representation
\begin{align}
\sigma x = Ax + Bu+E, \quad y=Cx+Du+F, \label{eq:state-space-affine} 
\end{align}
with 
$\scalebox{0.8}{$\bma\!
\begin{array}{ccc}
A & B & E  \\
C & D & F
\end{array}\!\ema $}
\in \R^{(n+p+1)\times (n+m+1)}$, using the same definitions given above for order and lag.

\begin{lemma}~\cite[Theorem 1]{martinelli2022data} \label{lemma:fundamental_generalized_affine}  
Let ${\B\in\A^q}$ be an ATI  system   of order $n$ and lag $\ell$ with state-space representation~\eqref{eq:state-space-affine}. Let ${w_d  = (u_d,y_d)   \in \R^{qT}}$  be a trajectory of  the system   and let ${L\in\mathbf{T}}$, with ${L > \ell}$. 
Then
\beq \label{eq:BL_affine} 
\B|_{[0,L-1]} = \aff \, H_L(w_d)
\eeq
if and only if
\begin{align}
\rank \,
	\bma 
		\begin{array}{c}  
		H_L(w_d) \\ 
		\mathbb{1}^{\transpose} 
		\end{array} 
	\ema  
&=  mL+n+1. \label{eq:generalized_persistence_of_excitation_affine}  
\end{align}
\end{lemma}

\begin{table*}[t]%
\captionsetup{width=0.9\hsize}
\bgroup
\renewcommand{\tabcolsep}{2.5pt}
\begin{tabular}{lclll} \hline
State-space representation  & Positive?  & Rank condition & Data-driven representation & Property  \\[0.3em] \hline
\rowcolor{Gray}
\scalebox{0.85}{$\sigma x = Ax + Bu, \, ~~~~~~~ y=Cx+Du$}
	& \newcrossmark
	& \scalebox{0.85}{$\rank \, H_L(w_d) =  mL+ n$}
	& \scalebox{0.85}{$\B|_{[0,L-1]} = \{ w  :  w = H_L(w_d) g, \, g\in\R^{T-L+1} \}  $}
 	& Linear  \\[0.3em] 
\scalebox{0.85}{$\sigma x = Ax + Bu, \, ~~~~~~~ y=Cx+Du$}
	& \newcheckmark  
	& \scalebox{0.85}{$\rank_{+} \, H_L(w_d) =  mL+ n$}
	& \scalebox{0.85}{$\B|_{[0,L-1]} = \{ w  : w = H_L(w_d)g , \, g\in\R_+^{T-L+1} \}  $}
 	& Polyhedral \\[0.3em] 
\rowcolor{Gray}
\scalebox{0.85}{$\sigma x = Ax + Bu + E, \, ~ y=Cx+Du+F$}
	& \newcrossmark  
 	& \scalebox{0.85}{$
\rank \,
	\bma 
		\begin{array}{c}  
		H_L(w_d) \\ 
		\mathbb{1}^{\transpose} 
		\end{array} 
	\ema  
=  mL+ n+1 $}
	& \scalebox{0.85}{$\B|_{[0,L-1]} = \{ w  :  w = H_L(w_d)g , \, g\in\R^{T-L+1}, \, \mathbb{1}^{\transpose} g = 1 \} $}
 	& Affine \\[0.3em] 
\scalebox{0.85}{$\sigma x = Ax + Bu + E, \, ~ y=Cx+Du+F$}
	& \newcheckmark  
 	& \scalebox{0.85}{$
\rank_+ \,
	\bma 
		\begin{array}{c}  
		H_L(w_d) \\ 
		\mathbb{1}^{\transpose} 
		\end{array} 
	\ema  
=  mL+ n+1 $}
	& \scalebox{0.85}{$\B|_{[0,L-1]} = \{ w  : w = H_L(w_d)g , \, g\in\R_+^{T-L+1}, \, \mathbb{1}^{\transpose} g = 1 \}$}
 	& Simplex \\[0.3em] \hline
\end{tabular}
\egroup
\vspace{0.1cm}
\caption{State-space representations of linear and affine systems with associated positivity constraints, generalized persistence of excitation rank conditions,  data-driven representations, and set-theoretic properties.}
\label{tab:data-driven-representations}
\end{table*}%

\noindent
To close this circle of ideas, we now establish analogous results for  \textit{positive}, linear and affine, time-invariant  systems with state-space representations~\eqref{eq:state-space-linear} and~\eqref{eq:state-space-affine}, respectively.

The analogy with Lemma~\ref{lemma:fundamental_generalized} and Lemma~\ref{lemma:fundamental_generalized_affine} is again surprisingly  simple. First,  one needs   to replace the usual notion of rank with that of  \textit{non-negative rank} of a matrix ${M\in\R_+^{p\times m}}$~\cite{cohen1993nonnegative}, defined as the smallest ${r\in\N}$ such that ${M=PQ}$, with ${P\in\R_+^{p\times r}}$ and ${Q\in\R_+^{r\times m}}$, and denoted by $\rank_+(M)$. 
Second, one needs to require that the state trajectory ${x_d \in \R^{nT}}$ corresponding to the input trajectory ${u_d \in \R^{mT}}$ is such that the matrix \scalebox{0.7}{$
\bma  
\begin{array}{c}
H_L(u_d) \\
H_1(x_d)
\end{array}
\ema  $}
has a monomial submatrix of order $mL+n$, i.e., $\rank_+\, \scalebox{0.7}{$
\bma  
\begin{array}{c}
H_L(u_d) \\
H_1(x_d)
\end{array}
\ema  $} = mL+n $ and all rows and columns of $\scalebox{0.7}{$
\bma  
\begin{array}{c}
H_L(u_d) \\
H_1(x_d)
\end{array}
\ema  $}$ have at most one non-zero entry~\cite[p.67]{berman1994nonnegative}.

\begin{theorem} \label{thm:fundamental_generalized_positive}  
Let ${\B\in\L^q}$ be an LTI,   positive  system of order $n$ and lag $\ell$ with state-space representation~\eqref{eq:state-space-linear}. 
Let ${w_d = (u_d,y_d)   \in \R^{qT}}$  be a trajectory of  the system, let ${x_d \in \R^{nT}}$ be the corresponding state trajectory, and let ${L\in\mathbf{T}}$, with ${L > \ell}$. 
Then
\beq 
\B|_{[0,L-1]} = \ccone \, H_L(w_d) 
\label{eq:BL_positive} 
\eeq
if and only  if 
\begin{align}
\rank_+ \,  H_L(w_d) =  mL+ n  \label{eq:generalized_persistence_of_excitation_positive}  
\end{align} 
and \scalebox{0.7}{$
\bma  
\begin{array}{c}
H_L(u_d) \\
H_1(x_d)
\end{array}
\ema  $}
has a monomial submatrix of order $mL+n$. 
\end{theorem}

\begin{theorem} \label{thm:fundamental_generalized_positive_affine}
 
Let ${\B\in\A^q}$ be an ATI,  positive system of order $n$ and lag $\ell$ with state-space  representation~\eqref{eq:state-space-affine}.    
Let ${w_d = (u_d,y_d)   \in \R^{qT}}$  be a trajectory of  the system, let ${x_d \in \R^{nT}}$ be the corresponding state trajectory, and let ${L\in\mathbf{T}}$, with ${L > \ell}$. 
Then
\beq \label{eq:BL_positive_affine} 
\B|_{[0,L-1]} = \conv \, H_L(w_d) 
\eeq
if and only if 
\beq
\rank_+ \,
	\bma 
		\begin{array}{c}  
		H_L(w_d) \\ 
		\mathbb{1}^{\transpose} 
		\end{array} 
	\ema  
=  mL+ n+1 \label{eq:generalized_persistence_of_excitation_affine_positive}  
\eeq 
and \scalebox{0.7}{$
\bma  
\begin{array}{c}
H_L(u_d) \\
H_1(x_d)
\end{array}
\ema  $}
has a monomial submatrix of order $mL+n$. 
\end{theorem}

\noindent 
Table~\ref{tab:data-driven-representations} summarizes the connections between the state-space representations~\eqref{eq:state-space-linear} and~\eqref{eq:state-space-affine}, positivity constraints, generalized persistence of excitation rank conditions,  data-driven representations, and set-theoretic properties.

Theorem~\ref{thm:fundamental_generalized_positive} and Theorem~\ref{thm:fundamental_generalized_positive_affine} have a number of practical implications. While computing the non-negative rank of a matrix is NP-hard~\cite{vavasis2010complexity}, there is a wealth of established algorithms that approximate the non-negative rank of a matrix within a certain factor of the true value~\cite{lee2000algorithms}. This can be useful in practical applications not only when an exact value is required, but also when an \textit{estimate} of the non-negative rank is sufficient (e.g., when the data is corrupted by noise).

Another important point is that the non-negative rank of a matrix is bounded below by its ordinary rank~\cite{cohen1993nonnegative}.  This is important in practice for input design or for situations where   
the non-negativity constrains can be disregarded.      In this case, one can generate non-negative persistently exciting signals (in an ordinary sense) and check \textit{a fortiori} that the data are informative, because the rank conditions~\eqref{eq:generalized_persistence_of_excitation} and~\eqref{eq:generalized_persistence_of_excitation_affine} imply
the non-negative rank conditions~\eqref{eq:generalized_persistence_of_excitation_positive} and~\eqref{eq:generalized_persistence_of_excitation_affine_positive}, respectively.

Another key point is the requirement for the matrix \scalebox{0.7}{$
\bma   \setlength{\arraycolsep}{2.5pt}
\begin{array}{c}
H_L(u_d) \\
H_1(x_d)
\end{array}
\ema  $}
to have a monomial submatrix of order ${mL+n}$. The condition may appear to be difficult to verify from data, since state trajectory measurements are typically not available. However, preliminary results indicate that this condition is guaranteed to hold for reachable positive systems~\cite{farina2000positive}, provided that $H_L(u_d)$ has a monomial matrix of order ${n+L}$.

Finally, we emphasize that persistence of excitation requirements imposed by non-negative rank conditions are \textit{weaker} than those imposed by standard rank conditions, despite being generally harder to compute. This is a reflection of the general principle that adding prior information about (the positivity of) a system facilitates the construction of models from data.
The next section provides an example to illustrate this point.

\section{Positivity-informed data-driven modeling} \label{sec:example}

Positivity of a system is crucial when evaluating if given data are informative for a prescribed model class.
We now illustrate this point by considering the Leslie model~\cite{leslie1945use}, a classical model of population ecology used to study the dynamics of population growth. The model is discrete-time and age-structured: the evolution of a population over time is described by dividing individuals into $n$     different age classes, indexed by ${i \in\mathbf{n}}$ and  ordered so that ${i = 1}$ is the class of newborns. The model relies on the assumption that the growth of the population depends on constant fertility rates, ${\alpha_i\in\Rge}$, and constant survival rates, $\beta_i \in [0,1]$, respectively. 
\bseq \label{eq:leslie}
The dynamics of the model are described by the equations
\begin{align}
\sigma x &= 
\scalebox{.85}{$\bma\!
\begin{array}{ccccc}
 \alpha_1 & \alpha_2 & \cdots & \alpha_{n-1} & \alpha_n \\
 \beta_1  & 0        & \cdots & 0 & 0 \\
 0        & \beta_2        & \ddots & \vdots & \vdots \\
 \vdots   & \ddots & \ddots & 0 & \vdots \\
 0        & \cdots & 0 & \beta_n & 0 
\end{array}
\ema$} x, 
\end{align}
where ${x_i(t)\in\Rge^n}$  is   the number of individuals in the age class ${i\in\mathbf{n}}$ 
at time ${t\in\Zge}$. We define the output as the total number of individuals in the last ${k\in\mathbf{n}}$ age classes, i.e., 
\beq
y = \, \scalebox{1}{$\bma\!
\begin{array}{ccccc}
 \,c_1 & \,c_2 & \ \cdots & \ \cdots \ & \, \ \, c_n
\end{array}
\ema$} x, 
\eeq
\eseq 
where ${c_{i}=0}$, if ${i\in\mathbf{n-k}}$, and ${c_i =1}$, otherwise.  The system is LTI, positive, and autonomous (i.e., ${m=0}$).  

For illustration, consider the Leslie model~\eqref{eq:leslie} defined by the matrices
\beq \nn
A = 
\bma
\begin{array}{cccc}
 0 & 0 & 0 & 1 \\
 1 & 0 & 0 & 0 \\
 0 & 1 & 0 & 0 \\
 0 & 0 & 1 & 0 
\end{array}
\ema , \quad
C = 
\bma
\begin{array}{cccc}
 0 & 0 & 1 & 1 
\end{array}
\ema .
\eeq
The order and the lag of the system are ${n=4}$ and ${\ell=4}$, respectively.
Now consider the trajectory ${w_d\in\Rge^7}$ of the system of length ${T=7}$ corresponding to 
the initial condition  ${x(0) = [\, 1 \ 0 \ 0 \ 0 \,]^{\transpose}}$, i.e.,
${w_d = [ \, 0 \ 0 \ 1 \ 1 \ 0 \ 0 \ 1  \, ]^{\transpose}}$.
The Hankel matrix of depth ${L=4}$ associated with the trajectory $w_d$ is
\beq \nn
H_4(w_d) = 
\bma
\begin{array}{cccc}
 0 & 0 & 1 & 1 \\
 0 & 1 & 1 & 0 \\
 1 & 1 & 0 & 0 \\
 1 & 0 & 0 & 1 
\end{array}
\ema
\eeq
The rank of the matrix $H_4(w_d)$ is
$\rank \, H_4(w_d) = 3$.
By Lemma~\ref{lemma:fundamental_generalized}, this implies that
the restricted behavior $\B|_{[0,3]}$ is not represented by  $\Span \, H_4(w_d)$.   
By contrast, the non-negative rank of the matrix $H_4(w_d)$ is
$\rank_+ \, H_4(w_d) =  4. $
This is because column permutations do not affect the non-negative rank of a matrix~\cite{dam1997unilaterally} and
because multiplying $H_4(w_d)$ by a suitable permutation matrix ${\Pi\in\R^{4\times 4}}$
yields 
\beq \nn
H_4(w_d) \Pi = 
\bma
\begin{array}{cccc}
 1 & 1 & 0 & 0 \\
 0 & 1 & 1 & 0 \\
 0 & 0 & 1 & 1 \\
 1 & 0 & 0 & 1 
\end{array} 
\ema ,
\eeq 
whose non-negative rank is  ${\rank_+ \, H_4(w_d) \Pi =  4}$~\cite[p.84]{berman1994nonnegative}. 
Furthermore, we have ${H_1(x_d) = I}$ and, hence, ${H_4(w_d) = \mathsf{O}_4}$ .  
By Theorem~\ref{thm:fundamental_generalized_positive}, we conclude that 
$\B|_{[0,3]} = \ccone \, H_4(w_d)$. As anticipated, this reflects the general principle that having prior information about (the positivity of) a system facilitates the construction of models from data.

\section{Conclusion} \label{sec:conclusion}

The paper has explored the role of conical and convex models in behavioral systems theory. 
Our analysis of their basic properties, as well as their construction from measured data, suggests 
that  imposing a conical and convex structure on a behavior leads to tractable representations, which retain many of the desirable properties of LTI models.

\appendix 
\section{Appendix}

\begin{proof}[Proof of Lemma~\ref{lemma:restrictions_and_time_shifts}]  
The first claim follows directly  from Lemma~\ref{lemma:algebraic_properties}, since the $t$-shift $\sigma^{t}$ is a linear map for any ${t \in \T}$.  For the second claim, assume $\B$ is a cone.  Fix ${\alpha \in \Rge}$ and ${t_1, t_2 \in \T},$ with ${t_1\le t_2}$. Then
\begin{align*}
{w \in \B|_{[t_1,t_2]}}
	&\,\Rightarrow \,{\exists\, v \in\B \,:\, w|_{[t_1,t_2]} = v|_{[t_1,t_2]}} \tag{def}\\
	&\,\Rightarrow \,{\alpha v \in\B} \tag{conicity}\\
	&\,\Rightarrow \,{\alpha w|_{[t_1,t_2]} = \alpha v|_{[t_1,t_2]} \in\B|_{[t_1,t_2]}} \tag{def}
\end{align*}
which shows that  $\B|_{[t_1,t_2]}$ is a cone. 
The proof is similar for convex (affine, linear) systems and it is thus omitted.
\end{proof}

\begin{proof}[Proof of Lemma~\ref{lemma:intersection_property}]
Let $\{\B_{j}\}_{j \in J}$ be a (possibly uncountably infinite) collection of 
closed, shift-invariant, conical behaviors.  
Then $\cap_{j \in J}  \B_j $ is clearly a cone and a closed set. 
Moreover, $\cap_{j \in J}  \B_j $ is shift-invariant, since
\beq \nn
\sigma^{t}\left( \cap_{j \in J}  \B_j \right)  \subseteq  
\cap_{j \in J}  \sigma^{t}\left( \B_j \right)  \subseteq 
\cap_{j \in J} \B_j, \text{  for all } {t \in \Zge}. 
\eeq
This proves that ${\cap_{j \in J} \B_j}$ is a closed shift-invariant cone and,
hence, that the model class $\K^q$ has the intersection property.
Similar arguments apply to the model classes $\C^q$ and $\A^q$.
\end{proof}

\begin{proof}[Proof of Theorem~\ref{thm:MPUM}]
By Lemma~\ref{lemma:intersection_property}, the model classes $\K^q$, $\C^q$, and $\A^q$ have the intersection property. 
By~\cite[Proposition 11]{willems1986timeb}, if a model class $\M$ has the intersection property and if 
${w_d\in(\R^q)^{\Zge}}$, the most powerful unfalsified model in the model class $\M$ based on the measurement ${w_d}$ is well-defined, unique, and given by ${\B_{\text{mpum}}(w_d) = \bigcap\limits_{\{\B \in \M: w_d \in \B \}} \B .}$
Thus, it suffices to verify this expression boils down to~\eqref{eq:MPUM_cone},~\eqref{eq:MPUM_conv}
and~\eqref{eq:MPUM_aff} for the model classes $\K^q$, $\C^q$, and $\A^q$, respectively.

Let ${\M = \K^q}$ and ${\tilde{\B} = \closure\left(\cone\{w_d, \sigma w_d, \sigma^2 w_d, \ldots \}\right)}$.
By definition, ${\tilde{\B}}$ is closed, shift-invariant, and conical, i.e., ${\tilde{\B} \in \K^q}.$ Furthermore, 
${w_d \in \tilde{\B}}.$ Finally,  if ${\B \in \K^q}$ is such that ${w_d \in \B}$,  then ${w_d \in \B}$ implies ${\cone\{w_d, \sigma w_d, \sigma^2 w_d, \ldots \} \subseteq \B}$ (by conicity and shift-invariance) and ${\closure(\cone\{w_d, \sigma w_d, \sigma^2 w_d, \ldots \}) \subseteq \B}$  (since $\B$ is closed),   that is, ${\tilde{\B} \subseteq \B}$.  
By uniqueness, we conclude that $\tilde{\B}$ is the most powerful unfalsified model in the model class $\K^q$ based on the measurement ${w_d}$. Similar considerations hold, \textit{mutatis mutandis}, for the model classes $\C^q$ and $\A^q$.
\end{proof}

\begin{proof}[Proof of Theorem~\ref{thm:fundamental_generalized_positive}]
Let ${L \in \mathbf{T}}$, with ${L>\ell}$. 
The proof relies on the key identity
\beq \label{eq:only_if_1}
\B|_{[0,L-1]} = \ccone \, ( M_L ).
\eeq 
where ${M_L \in \R^{(qL)\times (mL+n)}}$ is defined as
\beq \label{eq:M_L}
M_L = 
\Pi \,
\scalebox{0.95}{$\bma\!
\begin{array}{cc}
I            & 0 \\
\mathsf{T}_L & \mathsf{O}_L
\end{array}
\!\ema $},
\eeq
with ${\Pi\in \R^{(qL)\times (qL)}}$  a permutation matrix,
${\mathsf{O}_L\in\R^{(pL)\times n}}$ the observability matrix defined in~\eqref{eq:observability},   and 
${\mathsf{T}_L\in\R^{(pL)\times (mL)}}$ the (block-Toeplitz) \textit{convolution matrix}
\beq \label{eq:Toeplitz}
\mathsf{T}_L
= 
\scalebox{0.85}{$
\bma
\begin{array}{ccccc}
D  & 0 & \cdots  & 0 \\
CB & D & \ddots  & \vdots \\
\vdots & \ddots & \ddots  & 0 \\
C{A}^{L-2}B & \cdots & CB  & D
\end{array}
\ema $} .  
\eeq 
\noindent
Furthermore, we have
\beq \label{eq:only_if_7}  
H_L(w_d) = M_L 
\scalebox{0.8}{$
\bma  
\begin{array}{c}
H_L(u_d) \\
H_1(x_d)
\end{array}
\ema  $} ,
\eeq 
since $w_d$ is (by assumption) a trajectory of the system.

(Only if). Assume ${\B|_{[0,L-1]} = \ccone \, H_L(w_d).}$ Then
\beq \label{eq:only_if_11}
\ccone \, (M_L) 
\stackrel{\eqref{eq:only_if_1}}{=}   
\ccone \, (H_L(w_d))
\eeq
which, in view of~\cite[Theorem 14.1]{rockafellar1970convex}, is equivalent to 
\beq \label{eq:only_if_2}
\ccone \, (M_L)^{*} 
=
\ccone \, (H_L(w_d))^{*} ,
\eeq
where $C^{*}$ denotes the polar cone of the set $C$~\cite[p.121]{rockafellar1970convex}.
On the one hand, \eqref{eq:only_if_2} implies 
\beq \label{eq:only_if_3}
\ccone \, (M_L)^{*} \subseteq \ccone \, (H_L(w_d))^{*}.
\eeq
By Haar's Lemma~\cite[p.28]{dam1997unilaterally},~\eqref{eq:only_if_3} implies
\beq  \label{eq:only_if_HL}
H_L(w_d)  = M_L P  ,
\eeq
for some ${ P  \in\Rge^{(mL+n) \times (T-L+1)}}$. This, in turn, implies
\beq \label{eq:only_if_rank1}
\rank_+ \, H_L(w_d) \le mL+n.
\eeq 
On the other hand, \eqref{eq:only_if_2} implies 
\beq \label{eq:only_if_4}
\ccone \, (H_L(w_d))^{*} \subseteq \ccone \, (M_L)^{*} .
\eeq
By Haar's Lemma~\cite[p.28]{dam1997unilaterally},~\eqref{eq:only_if_4} implies
\beq  \label{eq:only_if_ML}
M_L = H_L(w_d) Q,
\eeq
for some ${ Q  \in\Rge^{(T-L+1)\times (mL+n)}}$.  Then \noindent
\bseq \label{eq:only_if_2_tot}
\begin{align}
mL+n 
& = \rank \, M_L  \label{eq:only_if_2a}  \\
& \le  \rank_+ \, M_L \label{eq:only_if_2b} \\
& \le \min\{\rank_+ \, H_L(w_d), \rank_+ \,  Q  \}  
\label{eq:only_if_2c} \\
& \le \min\{\rank_+ \, H_L(w_d), mL+n\},  \label{eq:only_if_2d}  
\end{align}
\eseq
where~\eqref{eq:only_if_2a} follows  from ${L>\ell}$,~\eqref{eq:only_if_2b} 
follows from~\cite[Lemma 2.3]{cohen1993nonnegative},~\eqref{eq:only_if_2c} follows from~\cite[Lemma 2.6]{cohen1993nonnegative}, and~\eqref{eq:only_if_2d} follows from $\min\{T-L+1, mL+n\} \le mL+n$. Then~\eqref{eq:only_if_rank1} and~\eqref{eq:only_if_2_tot} imply
\beq \label{eq:only_if_6}
\rank_+ \, H_L(w_d) = mL+n.
\eeq
Next, we show that the matrix 
$\scalebox{0.7}{$
\bma  \setlength{\arraycolsep}{2.5pt}
\begin{array}{c}
H_L(u_d) \\
H_1(x_d)
\end{array}
\ema  $}$
has a monomial submatrix of order $mL+n$. 
First, note that
\beq \nn
M_L 
\stackrel{\eqref{eq:only_if_ML}}{=} 
H_L(w_d) Q
\stackrel{\eqref{eq:only_if_7}  }{=} 
M_L 
\scalebox{0.8}{$
\bma  \setlength{\arraycolsep}{2.5pt}
\begin{array}{c}
H_L(u_d) \\
H_1(x_d)
\end{array}
\ema  $} Q,
\eeq
which implies
$
M_L (
\scalebox{0.8}{$
\bma   \setlength{\arraycolsep}{2.5pt}
\begin{array}{c}
H_L(u_d) \\
H_1(x_d)
\end{array}
\ema  $} Q -I)=0 .$
However, $\ker M_L=\{ 0\}$, since ${\rank \, M_L = mL+n}$. Then 
\beq \label{eq:only_if_Q_inverse}
\scalebox{0.8}{$
\bma   \setlength{\arraycolsep}{2.5pt}
\begin{array}{c}
H_L(u_d) \\
H_1(x_d)
\end{array}
\ema  $} Q = I .
\eeq
Thus, the full rank non-negative matrix $\scalebox{0.7}{$
\bma   \setlength{\arraycolsep}{2.5pt}
\begin{array}{c}
H_L(u_d) \\
H_1(x_d)
\end{array}
\ema  $}$ admits non-negative right inverse $Q$. By~\cite[Lemma 2]{plemmons1973regular}, this implies that 
$\scalebox{0.7}{$
\bma   \setlength{\arraycolsep}{2.5pt}
\begin{array}{c}
H_L(u_d) \\
H_1(x_d)
\end{array}
\ema  $}$
has a monomial submatrix of order $mL+n$.

(If)  Assume ${\rank_+ \, H_L(w_d) = mL+n}$
and  $\scalebox{0.7}{$
\bma   \setlength{\arraycolsep}{2.5pt}
\begin{array}{c}
H_L(u_d) \\
H_1(x_d)
\end{array}
\ema  $}$
has a monomial submatrix of order $mL+n$.
By Haar's Lemma~\cite[p.28]{dam1997unilaterally},~\eqref{eq:only_if_7} directly implies~\eqref{eq:only_if_3}. 
Furthermore, in view of ${\rank_+ \, M_L = mL+n}$ and~\cite[Lemma 2]{plemmons1973regular},
the (full rank) non-negative matrix $\scalebox{0.7}{$
\bma   \setlength{\arraycolsep}{2.5pt}
\begin{array}{c}
H_L(u_d) \\
H_1(x_d)
\end{array}
\ema  $}$ admits non-negative right inverse $Q$ such that~\eqref{eq:only_if_Q_inverse} holds.
Then \eqref{eq:only_if_7} and~\eqref{eq:only_if_Q_inverse} together imply~\eqref{eq:only_if_ML}. By Haar's Lemma~\cite[p.28]{dam1997unilaterally}, ~\eqref{eq:only_if_ML} implies~\eqref{eq:only_if_4}. By combining~\eqref{eq:only_if_3} and~\eqref{eq:only_if_4}, we obtain~\eqref{eq:only_if_2}, which in view of~\cite[Theorem 14.1]{rockafellar1970convex}, is equivalent to~\eqref{eq:only_if_11}. Finally, using~\eqref{eq:only_if_1}, we conclude that ${\B|_{[0,L-1]} = \ccone \, H_L(w_d)}$.
\end{proof}

\noindent
The proof of Theorem~\ref{thm:fundamental_generalized_positive_affine} is similar to that of Theorem~\ref{thm:fundamental_generalized_positive} and, hence, it is omitted in the interest of space.

\bibliographystyle{IEEEtran}
\bibliography{refs}

\end{document}